\documentclass[12pt]{article}

\usepackage{amsfonts,amssymb,amsmath,amsthm,epsfig,euscript}
\newtheorem{theorem}{Theorem}

\newtheorem{example}[theorem]{Example}

\newtheorem{conjecture}[theorem]{Conjecture}
\newtheorem{remark}[theorem]{Remark}

\long\def\symbolfootnote[#1]#2{\begingroup
\def\thefootnote{\fnsymbol{footnote}}\footnote[#1]{#2}\endgroup}
\def\S{\mathcal{S}}

\newcommand{\des}{\mathrm{des}}
\newcommand{\adj}{\mathrm{adj}}
\newcommand{\val}{\mathrm{val}}
\newcommand{\exc}{\mathrm{exc}}

\newcommand{\sg}{\sigma}

\def\E{\mathbb{E}}
\def\O{\mathbb{O}}
\def\N{\mathbb{N}}

\title{Equidistribution of $(X,Y)$-descents, $(X,Y)$-adjacent pairs, and $(X,Y)$-place-value pairs on permutations}

\author{Emeric Deutsch \\
\small Department of Mathematics\\[-0.8ex]
\small Polytechnic University\\[-0.8ex] \small Brooklyn,
NY 11201, USA. \\[-0.8ex]
\small \texttt{deutsch@duke.poly.edu} \and
Sergey Kitaev\footnote{The work presented here was supported by grant no. 090038011 from the Icelandic Research Fund.} \\
\small The Mathematics Institute\\[-0.8ex]
\small Reykjav\'{i}k University \\[-0.8ex]
\small IS-103 Reykjav\'{i}k, Iceland\\[-0.8ex]
\small \texttt{sergey@ru.is} \and
Jeffrey Remmel \\
\small Department of Mathematics\\[-0.8ex]
\small University of California, San Diego\\[-0.8ex]
\small La Jolla, CA 92093-0112. USA\\[-0.8ex]
\small \texttt{remmel@math.ucsd.edu} }

\date{\small Submitted: Date 1;  Accepted: Date 2;
 Published: Date 3.\\
\small MR Subject Classifications: 05A15}

\begin{document}

\maketitle

\begin{abstract}
An $(X,Y)$-descent in a permutation is a pair of adjacent elements
such that the first element is from $X$, the second element is from
$Y$, and the first element is greater than the second one. An
$(X,Y)$-adjacency in a permutation is a pair of adjacent elements
such that the first one is from $X$ and the second one is from $Y$.
An $(X,Y)$-place-value pair in a permutation is an element $y$ in
position $x$, such that $y$ is in $Y$ and $x$ is in $X$. It turns
out, that for certain choices of $X$ and $Y$ some of the three
statistics above become equidistributed. Moreover, it is easy to
derive the distribution formula for $(X,Y)$-place-value pairs thus
providing distribution for other statistics under consideration too.
This generalizes some results in the literature. As a result of our
considerations, we get combinatorial proofs of several remarkable
identities. We also conjecture existence of a bijection between two
objects in question preserving a certain statistic.
\end{abstract}

\section{Introduction}

Let $\S_n$ denote the set of permutations of $[n]=\{1, \ldots, n\}$
and $\N=\{1,2,\ldots\}$. Also, $\E$ and $\O$ denote the set of even
and odd numbers, respectively. For $\sg = \sg_1 \ldots \sg_n \in
\S_n$ define the following permutation statistics
\begin{eqnarray}
\des_{X,Y}(\sg) &=& |\{i: \sg_i>\sg_{i+1},\ \& \ \sg_i \in X \ \& \
\sg_{i+1} \in Y \}|,\nonumber \\ \adj_{X,Y}(\sg) &=& |\{i: \sg_i \in
X \ \& \ \sg_{i+1} \in Y \}|,\nonumber \\ \val_{X,Y}(\sg) &=& |\{i:
i \in X \ \& \ \sg_{i} \in Y \}|,\nonumber \\ \exc_{X,Y}(\sg) &=&
|\{i: \sg_i>i\ \&\ i\in X \ \& \ \sg_{i} \in Y
\}|,\nonumber\end{eqnarray} and the following corresponding
polynomials
\begin{eqnarray}
D_n^{X,Y}(x) &=& \sum_{\sg \in \S_n}
x^{\des_{X,Y}(\sg)}=\sum_{s=0}^{n-1} D_{n,s}^{X,Y}x^s, \nonumber\\
A_n^{X,Y}(x) &=& \sum_{\sg \in \S_n}
x^{\adj_{X,Y}(\sg)}=\sum_{s=0}^{n-1} A_{n,s}^{X,Y}x^s, \nonumber\\
V_n^{X,Y}(x) &=& \sum_{\sg \in \S_n}
x^{\val_{X,Y}(\sg)}=\sum_{s=0}^{n-1} V_{n,s}^{X,Y}x^s, \nonumber\\
E_n^{X,Y}(x) &=& \sum_{\sg \in \S_n}
x^{\exc_{X,Y}(\sg)}=\sum_{s=0}^{n-1} E_{n,s}^{X,Y}x^s. \nonumber
\end{eqnarray}
Objects counted by $\des_{X,Y}$ are called {\em $(X,Y)$-descents}
in~\cite{HR}. Similarly, we can talk of {\em $(X,Y)$-adjacencies},
{\em $(X,Y)$-place-value pairs}, and {\em $(X,Y)$-excedances}.

{\em Foata's first transformation}~\cite{F} exchanging excedances
and descents (to be used in the paper) can most easily be explained
with an example. The permutation $w=61437258$ has three excedances:
6, 4, and 7 in positions 1, 3, and 5, respectively. We write $w$ in
cycle form: (162)(34)(57)(8). Next, write each cycle with largest
element last, and order the cycles by increasing largest element:
(34)(216)(57)(8). Finally, reverse each cycle and erase the
parentheses to get the outcome permutation 43612758 with the
descents 43, 61, and 75.

\begin{remark}
Using Foata's first transformation, one obtains that
$D_n^{X,Y}(x)=E_n^{Y,X}(x)$. Thus, we do not need to provide any
arguments for the polynomial $E_n^{X,Y}(x)$ and its coefficients,
instead studying the other three polynomials.
\end{remark}

In this paper, we use the following notation for any $X
\subseteq \mathbb{N}$ and integer $ n\geq 1$:
\begin{eqnarray*}
X_n &=& [n] \cap X, x_n = |X_n|,
X_n^c = [n] -X,\ \mbox{and}\ x_n^c = |X_n^c|.\\
\end{eqnarray*}
Collecting some data on the polynomials, we noticed several
equidistributions among the statistics, and nice formulas associated
with them, for particular choices of sets $X$ and $Y$. We collect
those observations in Table~\ref{data}.
\begin{table}[htbp]
  \begin{center}\small{
      \begin{tabular}{|c|l|}
      \hline
 {\bf Stat.}& {\bf Description, related polynomial, and enumeration}\\
      \hline
$S_1$ & \# of even descent-tops ($D_{n,k}^{\E,\N}$). E.g., $S_1({\bf 2}15{\bf 4}36)=2$.\\
$S_2$ & \# of even excedance values ($E_{n,k}^{\N,\E}$). E.g., $S_2({\bf 2}15436)=1$.\\
$S_3$ & \# of even entries in even positions ($V_{n,k}^{\E,\E}$). E.g., $S_3(215{\bf 4}3{\bf 6})=2$.\\
& $a_{2n,k}=\left[n!{n \choose k}\right]^2$;
$a_{2n+1,k}=n!(n+1)!{n \choose k}{n+1 \choose k+1}$.\\
      \hline
$S_4$ & \# of odd descent-bottoms ($D_{n,k}^{\N,\O}$). E.g., $S_4(2{\bf 1}54{\bf 3}6)=2$.\\
$S_5$ & \# of odd excedance positions ($E_{n,k}^{\O,\N}$). E.g., $S_5({\bf 2}1{\bf 5}436)=2$.\\
$S_6$ & \# of even entries in odd positions ($V_{n,k}^{\O,\E}$). E.g., $S_6({\bf 2}15436)=1$.\\
$S_7$ & \# of odd entries in even positions ($V_{n,k}^{\E,\O}$). E.g., $S_7(2{\bf 1}5436)=1$.\\
$S_8$ & \# of (odd,even) pairs ($A_{n,k}^{\O,\E}$). E.g., $S_8(21\underline{54}\ \underline{36})=2$.\\
$S_9$ & \# of (even, odd) pairs ($A_{n,k}^{\E,\O}$). E.g.,
$S_9(\underline{21}5\underline{43}6)=2$.\\

& $a_{2n,k}=\left[n!{n \choose k}\right]^2$;
$a_{2n+1,k}=n!(n+1)!{n \choose k}{n+1 \choose k}$.\\
      \hline
$S_{10}$ & \# of odd descent-tops ($D_{n,k}^{\O,\N}$). E.g., $S_{10}(21{\bf 5}436)=1$.\\
$S_{11}$ & \# of odd excedance values ($E_{n,k}^{\N,\O}$). E.g., $S_{11}(21{\bf 5}436)=1$.\\
$S_{12}$ & \# of (odd,odd) pairs ($A_{n,k}^{\O,\O}$). E.g., $S_{12}(2\underline{15}436)=1$.\\
& $a_{2n,k}=(n!)^2{n-1 \choose k}{n+1\choose k+1}$;
$a_{2n+1,k}=n!(n+1)!{n \choose k}{n+1 \choose k}$.\\

      \hline
      $S_{13}$ & \# of even descent-bottoms ($D_{n,k}^{\N,\E}$). E.g., $S_{13}(215\bf{4}36)=1$.\\
$S_{14}$ & \# of even excedance positions ($E_{n,k}^{\E,\N}$). E.g., $S_{14}(215436)=0$.\\
& $a_{2n,k}=(n!)^2{n-1 \choose k}{n+1\choose k+1}$;
$a_{2n+1,k}=n!(n+1)!{n \choose k}{n+1 \choose k+1}$.\\

      \hline
$S_{15}$ & \# of odd entries in odd positions ($V_{n,k}^{\O,\O}$). E.g., $S_{15}(21{\bf 5}4{\bf 3}6)=2$.\\
& $a_{2n,k}=\left[n!{n \choose k}\right]^2$;
$a_{2n+1,k}=n!(n+1)!{n \choose k-1}{n+1 \choose k}$.\\
      \hline
$S_{16}$ & \# of (even,even) pairs ($A_{n,k}^{\E,\E}$). E.g., $S_{16}(215436)=0$.\\
& $a_{2n,k}=(n!)^2{n-1 \choose k}{n+1\choose k+1}$;
$a_{2n+1,k}=n!(n+1)!{n-1 \choose k}{n+2 \choose k+2}$.\\

      \hline
\end{tabular}}
  \end{center}
  \caption{16 statistics under consideration classified into 6 statistic groups.} \label{data}
\end{table}

Many of formulas listed in Table~\ref{data} are known (see,
e.g.,~\cite{KR}). Others are new but quite easy to prove. Our idea
to establish the equidistribution results is to prove general
recurrence relations for the statistics for arbitrary choice of sets
$X$ and $Y$. Then we will get the equidistributions in
Table~\ref{data} as a simple corollary to the fact that the
recurrences for the statistics in a given block are the same for a
particular choice of $X$ and $Y$. For example, we will show that
whenever $X$ and $Y$ are disjoint subsets of $\N$, then
$A_{n,s}^{X,Y} = V_{n,s}^{X,Y}$ for all $n$ and $s$. Indeed, the
recursions that we develop will allow us to give a bijective proof
of this fact. Other equidistribution results follow from simple
bijections. For example, it is easy to see that for any $X$ and $Y$,
$A_{n,s}^{X,Y} = A_{n,s}^{Y,X}$ since if $\sg_i \sg_{i+1}$ is an
$(X,Y)$-adjacency in $\sg = \sg_1 \ldots \sg_n$, then
$\sg_{i+1}\sg_i$ is a $(Y,X)$-adjacency in the reverse of $\sg$,
$\sg^r = \sg_n \sg_{n-1} \ldots \sg_1$.

Several of our formulas are quite easy to prove for one of our three
statistics. For example, it is always easy to compute
$V_{n,k}^{X,Y}$.

\begin{theorem}\label{formula} For any $X,Y \subseteq \mathbb{N}$,
\begin{equation}\label{eq:VXYnk}
V_{n,k}^{X,Y} = k! (x_n-k)!(x_n^c)!\binom{x_n}{k} \binom{y_n}{k} \binom{y_n^c}{x_n-k}.
\end{equation}
\end{theorem}

\begin{proof} To count the number of permutations of length $n$ with $k$ occurrences of $\val_{X,Y}$, we can first pick $k$
positions from $X_n$ in ${x_n\choose k}$ ways for the places where
we will have values of $Y$ occurring in the places corresponding to
$X_n$. Then we pick $k$ values from $Y_n$ in ${y_n\choose k}$ ways,
and permute the values in $k!$ ways to arrange the $k$ occurrences
of values in $Y_n$ in the places in $X_n$. In the remaining $x_n-k$
places in $X_n$, we must choose values from $Y_n^c$. We thus have
$\binom{y_n^c}{x_n-k}$ ways to choose those values and $(x_n-k)!$
ways to rearrange them. Finally we have $x_n^c!$ ways to arrange the
elements in places outside of $X_n$.
\end{proof}

Similarly, it is easy to count $A_{n,s}^{X,X}$ for any set $X \subseteq \N$.
That is, we have the following theorem.
\begin{theorem}\label{AXX} Suppose that $X \subseteq \N$ and for
any $n \geq 1$, $x_n=|X \cap [n]|$ and $x_n^c = |[n]-X|$.
Then for all $n \geq 1$,
\begin{equation}\label{eq:AXX}
A_{n,s}^{X,X} = (x_n)! (x_n^c)!\binom{x_n-1}{s}
\binom{x_n^c+1}{x_n-s}.
\end{equation}
\end{theorem}
\begin{proof} Fix $n \geq 1$. First we pick a permutation
$\sg$ of $X \cap [n]$ and a permutation $\tau$ of $[n]-X$. Clearly,
we have $(x_n)!(x_n^c)!$ ways to pick $\sg$ and $\tau$. We are now
interested in finding the number of permutations of $\gamma$ of
$S_n$ such that $\gamma$ restricted to the elements in $X \cap [n]$
yields the permutation $\sg$, $\gamma$ restricted to the elements in
$[n]-X$ yields the permutation $\tau$, and $\adj_{X,X}(\gamma) =s$.
Next in $\sg_1 \sg_2 \ldots \sg_{x_n}$, we think of choosing $s$
spaces from the $x_n-1$ spaces between the elements of $\sg$ to
create the adjacencies that will appear in such a $\gamma$.  For
example, if $n = 12$, $s=2$, $X = \E$, $\sg =4~2~10~8~6~12$, and we
pick spaces 2 and 5, then our choice  partitions $\sg$ into four
blocks, $4$, $2-10$, $8$ and $6-12$. Our idea is to insert these
blocks into the spaces that either lie immediately before an element
of $\tau$ or immediately after the last element
 of $\tau$. We label these spaces from left to right.
For example, suppose $\tau = 5~1~7~9~3~11$ and we pick spaces $2$,
$4$, $5$, and $7$. Then we would insert the block $4$ immediately
before $1$, the block $2-10$ immediately
before 9, the block $8$ immediately before $3$, $6-12$ immediately after
11 to obtain the permutation
$$ 5~4~1~7~2~10~9~8~3~11~6~12.$$
Clearly there are $\binom{x_n-1}{s}$ ways to choose the spaces to
obtain our $s$ adjacencies. This will leave us with $x_n-s$ blocks.
Then there are $\binom{x_n^c+1}{x_n-s}$ to choose the spaces for
$\tau$ where we insert the blocks.
\end{proof}

Hall and Remmel~\cite{HR} gave
direct combinatorial proofs of a pair of formulas for
$D_{n,s}^{X,Y}$ which combined with our equidistribution results,
gives formulas for the other polynomials under consideration. We
state these results here together with an example of using them.

\begin{theorem}\label{HR1}
\begin{equation}\label{Ieq:combXY1}
D_{n,s}^{X,Y} = \left| X_n^c\right|!\sum_{r=0}^s (-1)^{s-r} \binom{
\left| X_n^c\right| +r}{r} \binom{n+1}{s-r} \prod\limits_{x \in X_n}
(1+r + \alpha_{X,n,x} + \beta_{Y,n,x}),
\end{equation}
\end{theorem}

\begin{theorem}\label{HR2}
\begin{equation}\label{Ieq:combXY2}
D_{n,s}^{X,Y} = \left| X_n^c\right| !\sum_{r=0}^{ \left|
X_n\right|-s} (-1)^{\left| X_n\right|-s-r} \binom{\left|
X_n^c\right|+r}{r} \binom{n+1}{ \left| X_n \right|-s-r}
\prod\limits_{x \in X_n} (r+\beta_{X,n,x}-\beta_{Y,n,x}),
\end{equation}
where for any set $A$ and any $j,1 \leq j \leq n$, we define
\begin{eqnarray*}
\alpha_{A,n,j} &=&  |A^c \cap \{j+1, j+2,\ldots, n\}| = |\{x: j < x \leq n \ \& \ x \notin A\}|,\mbox{ and} \\
\beta_{A,n,j} &=&  |A^c \cap \{1, 2, \ldots, j-1\}| = |\{x: 1\leq  x
< j \ \& \ x \notin A\}|.
\end{eqnarray*}
\end{theorem}

\begin{example}
Suppose $X=\{2,3,4,6,7,9\}, Y = \{1, 4 , 8\}$, and $n=6$. Thus $X_6
= \{2, 3, 4, 6 \},X_6^c = \{1, 5 \}, Y_6 = \{1, 4 \}, Y_6^c = \{ 2,
3, 5, 6 \}$, and we have the following table of values of
$\alpha_{X,6,x}, \beta_{Y,6,x}$, and $\beta_{X,6,x}$.
\begin{center}
\begin{tabular}{|c|c|c|c|c|}
\hline
$x$ & $2$ & $3$ & $4$ & $6$ \\
\hline
$\alpha_{X,6,x}$ & $1$ & $1$ & $1$ & $0$ \\
\hline
$\beta_{Y,6,x}$ & $0$ & $1$ & $2$ & $3$ \\
\hline
$\beta_{X,6,x}$ & $1$ & $1$ & $1$ & $2$ \\
\hline
\end{tabular}
\end{center}
Equation (\ref{Ieq:combXY1}) gives
\begin{eqnarray*}
D_{6,2}^{X,Y} & = & 2!\sum\limits_{r=0}^2 (-1)^{2-r} \binom{2+r}{r} \binom{7}{2-r} (2+r)(3+r)(4+r)(4+r) \\
 & = & 2\left(1 \cdot 21 \cdot 2 \cdot 3 \cdot 4 \cdot 4 - 3 \cdot 7 \cdot 3 \cdot 4 \cdot 5 \cdot 5 + 6
\cdot 1 \cdot 4 \cdot 5 \cdot 6 \cdot 6 \right)
\\
 & = & 2(2016 - 6300 + 4320) \\
 & = & 72,
\end{eqnarray*}
while (\ref{Ieq:combXY2}) gives
\begin{eqnarray*}
D_{6,2}^{X,Y} & = & 2!\sum\limits_{r=0}^2 (-1)^{2-r} \binom{2+r}{r} \binom{7}{2-r} (1+r)(0+r)(-1+r)(-1+r) \\
 & = & 2 \left( 1 \cdot 21 \cdot 1 \cdot 0 \cdot (-1) \cdot (-1) - 3 \cdot 7 \cdot 2 \cdot 1 \cdot 0 \cdot 0
 + 6 \cdot 1 \cdot 3 \cdot 2 \cdot 1 \cdot 1 \right) \\
 & = & 2(0 - 0 + 36) \\
 & = & 72.
\end{eqnarray*}
\end{example}

The paper is organized as follows. In Section~\ref{sec2} we find
general recurrence relations for $D_{n,k}^{X,Y}$, $A_{n,k}^{X,Y}$,
and $V_{n,k}^{X,Y}$, and use them to explain the facts in
Table~\ref{data}. In Section~\ref{appl} we generalize several of the
results that appear in Table~\ref{data}, and use this to obtain
combinatorial proofs of several remarkable identities. Finally, in
Section~\ref{beyond}, we discuss some directions for further research.

\section{Recurrence relations for $D_{n,k}^{X,Y}$, $A_{n,k}^{X,Y}$, and
$V_{n,k}^{X,Y}$}\label{sec2}

In this section, we derive recurrence relations for $D_{n,k}^{X,Y}$,
$A_{n,k}^{X,Y}$, and $V_{n,k}^{X,Y}$. We notice that the recurrences
we get for $A_{n,k}^{X,Y}$ and $V_{n,k}^{X,Y}$ are almost identical,
except for the case when the element $n+1\in X\cap Y$ --- the
recurrences differ by ``1+.'' However, assuming $X\cap Y=\emptyset$,
we do not have this case, leading, in particular, to the explanation
of all of the equidistributions in Table~\ref{data}, and to many
more results for other choices of $X$ and $Y$, $X\cap Y=\emptyset$.

Another thing to observe is that in the case of the same recurrence
relations, we naturally get bijective proofs for the corresponding
equidistributed statistics. Indeed, one can label positions in a
permutation, say from left to right, in which we insert the largest
element, $n+1$, or do the other insertion procedure (see
Subsection~\ref{V}); then, it is enough to match insertions in the
positions having the same labels. However, such straightforward
approach is not necessarily the best one, as labeling positions
differently, rather than just from left to right, one may preserve
extra statistics in bijections (see Section~\ref{beyond} for
conjectures, which should be possible to prove using our approach
with different labeling).

\subsection{Recurrences for $D_{n,k}^{X,Y}$}

A recursion for $D_{n,k}^{X,Y}$ is derived in~\cite{HR}:
$$D_{n+1,k}^{X,Y}=\left\{
\begin{array}{ll} (k+1)D_{n,k+1}^{X,Y}+(n+1-k)D_{n,k}^{X,Y} & \mbox{ if }n+1\not\in X, \\
(y_n-(k-1))D_{n,k-1}^{X,Y}+(n+1-(y_n-k))D_{n,k}^{X,Y} & \mbox{ if
}n+1\in X.
\end{array}\right.$$ An argument for deriving the recursion is as follows. We are thinking of inserting the element $n+1$ in a permutation
$\sg=\sg_1\ldots\sg_n$, and we consider which of the obtained
permutations are counted by $D_{n,k}^{X,Y}$. If $n+1\not\in X$ then
one never increases the number of $(X,Y)$-descents by inserting
$n+1$. More precisely, the number of $(X,Y)$-descents is either
unchanged, or it is decreased by 1, when $n+1$ is inserted between
$\sg_i\in X$ and $\sg_{i+1}\in Y$ where $\sg_i>\sg_{i+1}$. The
corresponding recursion case follows.

For the second case, notice that if $n+1\in X$, then the number of
$(X,Y)$-descents is unchanged if $n+1$ is inserted at the end of the
permutation, in front of $\sg_j\not\in Y$, or between $\sg_i\in X$
and $\sg_{i+1}\in Y$ where $\sg_i>\sg_{i+1}$, and it is increased by
1 in other cases (that is, when $n+1$ is inserted in front of
$\sg_j\in Y$ not involved in an $(X,Y)$-descent). The second
recursion case follows.

We use a similar approach to derive recurrence relations for
$A_{n,k}^{X,Y}$. Our derivations for $V_{n,k}^{X,Y}$ use a different
insertion procedure.

\subsection{Recurrences for $A_{n,k}^{X,Y}$}

We consider 4 cases.

{\bf Case 1.} $n+1\not\in X\cup Y$. The number of $(X,Y)$-adjacent
pairs is decreased by 1 when $n+1$ is inserted between $\sg_i\in X$
and $\sg_{i+1}\in Y$ and it is unchanged otherwise. Thus, in this
case
$$A_{n+1,k}^{X,Y}=(k+1)A_{n,k+1}^{X,Y}+(n+1-k)A_{n,k}^{X,Y}.$$

{\bf Case 2.} $n+1\in X\cap Y$. Adding $n+1$ after a $\sigma_i\in X$
or before a $\sigma_j\in Y$ increases $\adj_{X,Y}$ by 1, while it
keeps $\adj_{X,Y}(\sg)$ unchanged otherwise. However, we note that
the place between $\sg_i\in X$ and $\sg_{i+1}\in Y$ is after a
$\sg_i\in X$ {\em and} before a $\sg_{i+1}\in Y$. Thus, in this case
$$A_{n+1,k}^{X,Y}=(x_n+y_n-(k-1))A_{n,k-1}^{X,Y}+(n+1-(x_n+y_n-k))A_{n,k}^{X,Y}.$$

{\bf Case 3.}  $n+1\in X-Y$. Inserting $n+1$ to the left of a
$\sg_i\not\in Y$ does not change $\adj_{X,Y}(\sg)$, which is also
the case if $n+1$ is inserted between $\sg_i\in X$ and $\sg_{i+1}\in
Y$, or $n+1$ is inserted at the very end. On the other hand, if
$n+1$ is inserted between $\sg_i\not\in X$ and $\sg_{i+1}\in Y$, the
number of $(X,Y)$-adjacent pairs is increased by 1. Thus, in this
case
$$A_{n+1,k}^{X,Y}=(y_n-(k-1))A_{n,k-1}^{X,Y}+(n+1-(y_n-k))A_{n,k}^{X,Y}.$$

{\bf Case 4.}  $n+1\in Y-X$. Inserting $n+1$ to the right of a
$\sg_i\not\in X$ does not change $\adj_{X,Y}(\sg)$, which is also
the case if $n+1$ is inserted between $\sg_i\in X$ and $\sg_{i+1}\in
Y$, or $n+1$ is inserted at the very beginning. On the other hand,
if $n+1$ is inserted between $\sg_i\in X$ and $\sg_{i+1}\not\in Y$,
the number of $(X,Y)$-adjacent pairs is increased by 1. Thus, in
this case
$$A_{n+1,k}^{X,Y}=(x_n-(k-1))A_{n,k-1}^{X,Y}+(n+1-(x_n-k))A_{n,k}^{X,Y}.$$

\subsection{Recurrences for $V_{n,k}^{X,Y}$}\label{V}

Instead of inserting the largest element, $n+1$, in all possible
places, we use another insertion procedure $I_n^{(i)}(\sg)$ that
generates $\S_{n+1}$ from $\S_{n}$. For $\sg=\sg_1\ldots\sg_n$, let
$I_{n+1}^{(n+1)}(\sg)=\sg(n+1)=\sg_1\ldots\sg_n(n+1)$, and for
$1\leq i\leq n$, let
$I_{n+1}^{(i)}(\sg)=\sg_1\ldots\sg_{i-1}(n+1)\sg_{i+1}\ldots\sg_n\sg_i$
(that is, in the last case we replace $\sg_i$ in $\sg$ by $n+1$ and
move $\sg_i$ to the very end).

We now consider 4 cases.

{\bf Case 1.} $n+1\not\in X\cup Y$. In this case, one can only
decrease the number of $(X,Y)$-place-value pairs. This happens when
$n+1$ occupies position $i\in X$ in $I_{n+1}^{(i)}(\sg)$ for some
$\sg$, such that $\sg_i\in Y$ ($\sg_i$ is in position $n+1$ in
$I_{n+1}^{(i)}(\sg)$). Thus, in this case
$$V_{n+1,k}^{X,Y}=(k+1)V_{n,k+1}^{X,Y}+(n+1-k)V_{n,k}^{X,Y}.$$

{\bf Case 2.} $n+1\in X\cap Y$. This is straightforward to see that
the number of $(X,Y)$-place-value pairs is unchanged if $i\not\in X$
and $\sg_i\not\in Y$, and it increases by 1 in each of the following
three cases: $i\in X$ and $\sigma_i\in Y$, $i\in X$ and
$\sigma_i\not\in Y$, and $i\not\in X$ and $\sigma_i\in Y$. Note that
we add 1 for each $i\in X$ and 1 for each $\sg_i\in Y$, so we count
$i\in X$ and $\sg_i\in Y$ twice. Moreover, having $n+1$ in position
$n+1$ gives one more $(X,Y)$-place-value pair. Thus, in this case
$$V_{n+1,k}^{X,Y}=(1+x_n+y_n-(k-1))V_{n,k-1}^{X,Y}+(n+1-(1+(x_n+y_n-k)))V_{n,k}^{X,Y}.$$

{\bf Case 3.}  $n+1\in X-Y$. One can check that in this case, the
number of $(X,Y)$-place-value pairs increases by 1 if $i\not\in X$
and $\sigma_i\in Y$, and it is unchanged otherwise. Thus, in this
case
$$V_{n+1,k}^{X,Y}=(y_n-(k-1))V_{n,k-1}^{X,Y}+(n+1-(y_n-k))V_{n,k}^{X,Y}.$$

{\bf Case 4.}  $n+1\in Y-X$.  One can check that in this case, the
number of $(X,Y)$-place-value pairs increases by 1 if $i\in X$ and
$\sigma_i\not\in Y$, and it is unchanged otherwise. Thus, in this
case
$$V_{n+1,k}^{X,Y}=(x_n-(k-1))V_{n,k-1}^{X,Y}+(n+1-(x_n-k))V_{n,k}^{X,Y}.$$
\ \\
\ \\
There are a number of cases where the recursions for
$V_{n,k}^{A,B}$, $A_{n,k}^{C,D}$, and $D_{n,k}^{E,F}$ coincide so
that we immediately have equality between the various pairs of
statistics.  For example, comparing the
recursions for $A_{n,k}^{X,Y}$ and $V_{n,k}^{X,Y}$, we immediately
have the following theorem.

\begin{theorem}\label{V=A}
If $X$ and $Y$ are subsets of $\N$ such that $X \cap Y = \emptyset$,
then for all $n$ and $k$,
$V_{n,k}^{X,Y} = A_{n,k}^{X,Y}$.
\end{theorem}

In fact, it is easy to see that our
proofs of the recursions can be used to give inductive proof
that there exists a  bijection from $S_n$ onto $S_n$ for all $n$
that will witness this equality. That is,
our proofs of the recursions immediately allow us to construct
inductively bijections $\Theta_n:S_n \rightarrow S_n$ for all $n$
such that for all $\sg \in S_n$,
$$\adj_{X,Y}(\sg) = \val_{X,Y}(\Theta_n(\sg)).$$
For example, suppose that we have constructed $\Theta_n$ and $n+1
\not \in X \cup Y$. First consider our insertion procedure to prove
the recursions for $A_{n,s}^{X,Y}$. If $\sigma \in S_n$, then we
consider the places where we can insert $n+1$ to $\sigma$. We first
label the spaces between the elements $\sg_i \in X$ and $\sg_{i+1}
\in Y$ from left to right with $1, \ldots, \adj_{X,Y}(\sg)$ and then
label the rest of the spaces from left to right with
$\adj_{X,Y}(\sg)+1, \ldots , n+1$. For example, if $X = \E$, $Y =
\O$, and $\sg = 1~4~3~2~5$, the spaces would be labeled by
$$_{\stackrel{-}{3}}1_{\stackrel{-}{4}}4_{\stackrel{-}{1}}3_{\stackrel{-}{5}}
2_{\stackrel{-}{2}}5_{\stackrel{-}{6}}.$$ We then let $\sg^{(i)}$ be
the permutation that results by inserting $n+1$ into the space
labeled $i$. For example, in our example, $\sg^{(4)} = 1~6~4~3~2~5$.
Next we consider our insertion procedure for proving the recursions
for $V_{n,s}^{X,Y}$. Now if $\tau \in S_n$, then we label the
positions of $\tau$ by first labeling the positions $i$ such that $i
\in X$ and $\tau_i \in Y$ from left to right with $1, \ldots
,\val_{X,Y}(\tau)$ and then label the remaining positions from left
to right with $\val_{X,Y}(\tau)+1, \ldots, n.$ For example, if $X =
\O$ and $Y = \E$ and $\tau = 1~4~2~5~3$, then we would label the
positions
$$\frac{1}{{\bf 2}}~\frac{4}{{\bf 3}}~\frac{2}{{\bf 1}}~\frac{5}{{\bf 4}}~\frac{3}{{\bf 5}}$$
where we have indicated the labels in boldface.
If label $j$ is in position $i$, then we let
$\tau^{(j)} = I_n^{(i)}(\tau)$ and we let
$\tau^{(n+1)} = I_n^{(n+1)}(\tau)$. For example, in our
case, $\tau^{(2)} = 6~4~2~5~3~1$. Then for any $\sg \in S_n$ and
$i \in \{1, \ldots, n+1\}$,
we can define
$$\Theta_{n+1}(\sg^{(i)}) = \Theta_n(\sg)^{(i)}.$$
We can extend $\Theta_n$ to $\Theta_{n+1}$ in the other cases of the
recursions in a similar manner.

Similarly, comparing the recursions for the $V_{n,k}^{A,B}$,
$A_{n,k}^{C,D}$, and $D_{n,k}^{E,F}$, we can also derive bijective
proofs of the following theorems.

\begin{theorem}\label{D=A} If $X$ and $Y$ are subsets of $\N$,
$A= X \cup Y$ and there exists a $B \subseteq \N$ such that
$b_n =|B \cap [n]|$ satisfies
$$b_n = \begin{cases}
x_n+y_n =|X \cap [n]| + |Y \cap [n]| & \ \mbox{if} \ n+1 \in X \cap Y \\
y_n = |Y \cap [n]|& \ \mbox{if} \ n+1 \in X  - Y \\
x_n = |X \cap [n]|& \ \mbox{if} \ n+1 \in Y - X,
\end{cases}
$$
then $D_{n,k}^{A,B} = A_{n,k}^{X,Y}$.
\end{theorem}

\begin{theorem}\label{D=V} If $X$ and $Y$ are subsets of $\N$,
$A= X \cup Y$ and there exists a $B \subseteq \N$ such that
$b_n =|B \cap [n]|$ satisfies
$$b_n = \begin{cases}
1+x_n+y_n = 1+ |X \cap [n]|+ |Y \cap [n]|& \ \mbox{if} \ n+1 \in X \cap Y \\
y_n = |Y \cap [n]|& \ \mbox{if} \ n+1 \in X  - Y \\
x_n =|X \cap [n]| & \ \mbox{if} \ n+1 \in Y - X,
\end{cases}
$$
then $D_{n,k}^{A,B} = V_{n,k}^{X,Y}$.
\end{theorem}

\subsection{Explanation of Table~\ref{data} using our general results}

\begin{enumerate}
\item {\bf The first group of statistics.}
$D_{n,k}^{\E,\N}=E_{n,k}^{\N,\E}$ by Foata's first transformation.
Also, $D_{n,k}^{\E,\N}=V_{n,k}^{\E,\E}$ by Theorem~\ref{D=V}.
Indeed, in this case $A=X=Y=\E$ and $B=\N$ leading to $A=X\cup Y$,
$X-Y=Y-X=\emptyset$, and $b_n=n=1+2|\E\cap [n]|$ if $n+1\in \E$. As
for the formulas, we can apply Theorem~\ref{formula} with $X=Y=\E$:
$$a_{2n,k}=V_{2n,k}^{X,Y}=k!(n-k)!n!{n \choose k}{n\choose k}{n\choose n-k}=\left[n!{n \choose
k}\right]^2;$$
$$a_{2n+1,k}=V_{2n+1,k}^{X,Y}=k!(n-k)!(n+1)!{n\choose k}{n\choose k}{n+1\choose n-k}=n!(n+1)!{n \choose k}{n+1 \choose k+1}.$$

\item {\bf The second group.} $D_{n,k}^{\N,\O}=E_{n,k}^{\O,\N}$ by Foata's
first transformation. Applying the reverse operation to each
permutation, one sees that $A_{n,k}^{\O,\E}=A_{n,k}^{\E,\O}$.
Applying the inverse operation to each permutation, one gets
$V_{n,k}^{\O,\E}=V_{n,k}^{\E,\O}$. By Theorem~\ref{V=A},
$V_{n,k}^{\O,\E}=A_{n,k}^{\O,\E}$ as $\O\cap\E=\emptyset$. Finally,
by Theorem~\ref{D=A}, $D_{n,k}^{\N,\O}=A_{n,k}^{\O,\E}$. Indeed, in
this case $A=\N$, $B=X=\O$, and $Y=\E$ leading to $A=X\cup Y$, and
$$b_n=\mbox{\# of odd numbers in }[n]=\left\{
\begin{array}{ll} \E\cap [n] & \mbox{ if }n+1\not\in \O, \\
\O \cap [n] & \mbox{ if }n+1\in \E.
\end{array}\right.$$
As for the formulas, we can apply Theorem~\ref{formula} with $X=\E$
and $Y=\O$:
$$a_{2n,k}=V_{2n,k}^{\O,\E}=k!(n-k)!n!{n \choose k}{n\choose k}{n\choose n-k}=\left[n!{n \choose
k}\right]^2;$$
$$a_{2n+1,k}=V_{2n+1,k}^{\O,\E}=k!(n-k)!(n+1)!{n\choose k}{n+1\choose k}{n\choose n-k}=n!(n+1)!{n \choose k}{n+1 \choose k}.$$

\item {\bf The third group.} Again, $D_{n,k}^{\O,\N}=E_{n,k}^{\N,\O}$ by Foata's
first transformation. Moreover, by Theorem~\ref{D=A},
$D_{n,k}^{\O,\N}=A_{n,k}^{\O,\O}$. Indeed, in this case $A=X=Y=\O$
and $B=\N$ leading to $A=X\cup Y$, $X-Y=Y-X=\emptyset$, and
$b_n=n=2|\O\cap [n]|$ if $n+1\in \O$. As for the formulas, we can
apply Theorem~\ref{AXX} with $X=\O$:
$$a_{2n,k}=A_{2n,k}^{\O,\O}=(n!)^2{n-1 \choose k}{n+1\choose k+1};$$
$$a_{2n+1,k}=A_{2n+1,k}^{\O,\O}=n!(n+1)!{n \choose k}{n+1 \choose k}.$$

\item {\bf The fourth group.} $D_{n,k}^{\N,\E}=E_{n,k}^{\E,\N}$ by Foata's
first transformation. The formulas for $D_{n,k}^{\N,\E}$ are proved
in~\cite[Section 4]{KR}.
\item {\bf The fifth group.} We use Theorem~\ref{formula} with
$X=Y=\O$, to get
$$a_{2n,k}=V_{2n,k}^{\O,\O}=k!n!n!{n \choose k}{n\choose k}{n\choose n-k}=\left[n!{n \choose k}\right]^2;$$
$$a_{2n+1,k}=V_{2n+1,k}^{\O,\O}=k!(n+1-k)!n!{n+1\choose k}{n+1\choose k}{n\choose k-1}=n!(n+1)!{n \choose k-1}{n+1 \choose k}.$$
\item {\bf The sixth group.} We use Theorem~\ref{AXX} with
$X=\E$, to get
$$a_{2n,k}=A_{2n,k}^{\E,\E}=(n!)^2{n-1 \choose k}{n+1\choose k+1};$$
$$a_{2n+1,k}=A_{2n+1,k}^{\E,\E}=n!(n+1)!{n-1\choose k}{n+2\choose n-k}.$$

\end{enumerate}

\section{Applications}\label{appl}

In this section, we shall generalize several of the results that
appear in Table~\ref{data}. That is, in Table~\ref{data}, we
consider the parity of the elements in a descent, adjacency, or
place-value pair. We shall show that we can get similarly formulas
when we consider the equivalence class modulo $k$ of the elements in
a descent, adjacency, or place-value pair. See~\cite{KR0} for
related research on descents generalizing results of~\cite{KR}. For
any $k \geq 2$ and $0 \leq i \leq k-1$, we let $i+k\N = \{i+kn: n
\geq 0\}$.

First, we shall consider $V_{n,k}^{X,Y}$ and $A_{n,k}^{X,Y}$ where
$X = i+k\N$ and $Y = j+k\N$ and $0 \leq i < j \leq k-1$. It follows
from Theorem \ref{V=A} that $V_{n,k}^{X,Y} = A_{n,k}^{X,Y}$ in the
case. Suppose that $A =i+k\N \cup j+k\N$ and $B = i + k\N$. Note
that when $m+1 = kn+i \in X -Y$, then $y_m = n = b_m =|B\cap [m]|$
and when $m+1 = kn+j \in Y -X$, then $x_m = n+1 = b_m$. Thus it
follows from Theorems \ref{D=A}  and \ref{D=V} that $D_{n,k}^{A,B} =
V_{n,s}^{X,Y} = A_{n,s}^{X,Y}$ for all $n$ and $s$.
We then have three cases.\\
\ \\
{\bf Case 1.} $m = kn+t$ where $0 \leq t <i$. In this case, $x_m =|X
\cap [m]| = y_m = |Y \cap [m]| = n$ and
 $x_m^c =|[m]-X| = y_m^c = |[m]-Y| = (k-1)n+t$. Thus it follows
from Theorem \ref{formula} that
\begin{equation}\label{ap1}
V_{m,s}^{X,Y} = \binom{n}{s}^2 \binom{(k-1)n+t}{n-s} s! (n-s)!
((k-1)n+t)!.
\end{equation}
On the other hand, it follows from Theorem \ref{HR1} that
\begin{equation}\label{ap2}
D_{m,s}^{A,B} = |A_m^c|! \sum_{r=0}^s (-1)^{s-r} \binom{|A_m^c|+r}{r}
\binom{m+1}{s-r} \prod_{x \in A_m} (1+r +\alpha_{A,m,x} +\beta_{B,m,x}).
\end{equation}
In this case $|A_m^c| = kn+t-2n = (k-2)n+t$.
For any $x$, it is easy to see that
$$\alpha_{A,m,x} +\beta_{B,m,x} = kn+t-1  -|A \cap [x+1,kn+t] -|B\cap[x-1]|$$
where $[x+1,kn+t] = \{r:x+1 \leq r \leq kn+t\}$. Thus
for any $0 \leq a \leq n-1$,
\begin{eqnarray}\label{ap3}
\alpha_{A,m,ak+i} +\beta_{B,m,ak+i} &=&  kn+t-1  - (2n -(2a+1)) - a
\nonumber \\
&=& (k-2)n+t +a
\end{eqnarray}
and
\begin{eqnarray}\label{ap4}
\alpha_{A,m,ak+j} +\beta_{B,m,ak+j} &=&  kn+t-1  - (2n -(2a+2)) - (a+1)
\nonumber \\
&=& (k-2)n+t +a.
\end{eqnarray}
Thus
\begin{eqnarray*}
&&\prod_{x \in A_m} (1+r +\alpha_{A,m,x} +\beta_{B,m,x}) = \\
&& \left(\prod_{a=0}^{n-1} (1+r +\alpha_{A,m,ak+i}
+\beta_{B,m,ak+i})\right)
 \left(\prod_{a=0}^{n-1} (1+r +\alpha_{A,m,ak+j} +\beta_{B,m,ak+j})\right) =\\
&& \left(\prod_{a=0}^{n-1} (1+r+(k-2)n+t+a)\right)^2 = \\
&&(1+r+(k-2)n+t)_n (1+r+(k-2)n+t)_n
\end{eqnarray*}
where we define $(a)_n$ by $(a)_0 = 1$ and $(a)_n = a(a+1) \cdots
(a+n-1)$ for $n \geq 1$. Since we have a combinatorial proof of the
fact that $V_{m,s}^{X,Y} =D_{m,s}^{A,B}$ in this case and the proof
of Theorem \ref{HR1} is also completely combinatorial, it follows
that we have a combinatorial proof of the following identity:
\begin{eqnarray}
&&\binom{n}{s}^2 \binom{(k-1)n+t}{n-s} s! (n-s)!
((k-1)n+t)! = \\
&& ((k-2)n+t)! \sum_{r=0}^s (-1)^{s-r}
\binom{(k-2)n+t +r}{r} \binom{kn+t+1}{s-r} \times \nonumber \\
&& \ \ \ \ \ \ \ \ \ \ (1+r+(k-2)n+t)_n (1+r+(k-2)n+t)_n. \nonumber
\end{eqnarray}
\ \\
{\bf Case 2.} $m = kn+t$ where $i \leq t <j$. In this case, $x_m =|X
\cap [m]| = n+1$ and $y_m = |Y \cap [m]| = n$ and
 $x_m^c =|[m]-X| = (k-1)n+t-1$ and
$y_m^c = |[m]-Y| = (k-1)n+t$. Thus it follows
from Theorem \ref{formula} that
\begin{equation}\label{ap21}
V_{m,s}^{X,Y} = \binom{n+1}{s}\binom{n}{s} \binom{(k-1)n+t}{n+1-s} s! (n+1-s)!
((k-1)n+t-1)!.
\end{equation}

On the other hand, we can obtain a formula for $V_{m,s}^{X,Y}=
D_{m,s}^{A,B}$ from equation (\ref{ap2}). In this case $|A_m^c| =
kn+t-(2n+1) = (k-2)n+t-1$. For any $0 \leq a \leq n$,
\begin{eqnarray}\label{ap23}
\alpha_{A,m,ak+i} +\beta_{B,m,ak+i} &=&  kn+t-1  - (2n+1 -(2a+1)) - a
\nonumber
\\
&=& (k-2)n+t-1 +a
\end{eqnarray}
and, for any $0 \leq a \leq n-1$
\begin{eqnarray}\label{ap24}
\alpha_{A,m,ak+j} +\beta_{B,m,ak+j} &=&  kn+t-1  - (2n+1 -(2a+2)) - (a+1)
\nonumber \\
&=& (k-2)n+t-1 +a.
\end{eqnarray}
Thus
\begin{eqnarray*}
&&\prod_{x \in A_m} (1+r +\alpha_{A,m,x} +\beta_{B,m,x}) = \\
&& \left(\prod_{a=0}^{n} (1+r +\alpha_{A,m,ak+i}
+\beta_{B,m,ak+i})\right)
 \left(\prod_{a=0}^{n-1} (1+r +\alpha_{A,m,ak+j} +\beta_{B,m,ak+j})\right) =\\
&&(1+r+(k-2)n+t-1)_{n+1} (1+r+(k-2)n+t-1)_n =\\
&&(r+(k-2)n+t)_{n+1} (r+(k-2)n+t)_n.
\end{eqnarray*}
As in Case 1, it follows that we have a combinatorial proof of
the following identity:
\begin{eqnarray}
&&\binom{n+1}{s}\binom{n}{s} \binom{(k-1)n+t}{n+1-s} s! (n-s)!
((k-1)n+t-1)! = \\
&& ((k-2)n+t-1)! \sum_{r=0}^s (-1)^{s-r}
\binom{(k-2)n+t-1 +r}{r} \binom{kn+t+1}{s-r} \times \nonumber \\
&&\ \ \ \ \ \ \ \ \ \ (r+(k-2)n+t)_{n+1} (r+(k-2)n+t)_n. \nonumber
\end{eqnarray}
\ \\
{\bf Case 3.} $m = kn+t$ where $j \leq t \leq k-1$. In this case,
$x_m =|X \cap [m]| = y_m = |Y \cap [m]| = n+1$ and
 $x_m^c =|[m]-X| = y_m^c = |[m]-Y| = (k-1)n+t-1$. Thus it follows
from Theorem \ref{formula} that
\begin{equation}\label{ap31}
V_{m,s}^{X,Y} = \binom{n+1}{s}^2 \binom{(k-1)n+t-1}{n+1-s} s!
(n+1-s)! ((k-1)n+t-1)!.
\end{equation}

On the other hand, we can obtain a formula for $V_{m,s}^{X,Y}=
D_{m,s}^{A,B}$ from equation (\ref{ap2}). In this case $|A_m^c| =
kn+t-(2n+2) = (k-2)n+t-2$. For any $0 \leq a \leq n$,
\begin{eqnarray}\label{ap33}
\alpha_{A,m,ak+i} +\beta_{B,m,ak+i} &=&  kn+t-1  - (2n+2 -(2a+1)) - a
\nonumber \\
&=& (k-2)n+t-2 +a
\end{eqnarray}
and, for any $0 \leq a \leq n$
\begin{eqnarray}\label{ap34}
\alpha_{A,m,ak+j} +\beta_{B,m,ak+j} &=&  kn+t-1  - (2n+2 -(2a+2)) - (a+1)
\nonumber \\
&=& (k-2)n+t-2 +a.
\end{eqnarray}
Thus
\begin{eqnarray*}
&&\prod_{x \in A_m} (1+r +\alpha_{A,m,x} +\beta_{B,m,x}) = \\
&& \left(\prod_{a=0}^{n} (1+r +\alpha_{A,m,ak+i}
+\beta_{B,m,ak+i})\right)
 \left(\prod_{a=0}^{n} (1+r +\alpha_{A,m,ak+j} +\beta_{B,m,ak+j})\right) =\\
&&(1+r+(k-2)n+t-2)_{n+1} (1+r+(k-2)n+t-2)_{n+1} =\\
&&(r+(k-2)n+t-1)_{n+1} (r+(k-2)n+t-1)_{n+1}.
\end{eqnarray*}
It follows that we have a combinatorial proof of the following
identity.

\begin{eqnarray}
&&\binom{n+1}{s}^2 \binom{(k-1)n+t-1}{n+1-s} s! (n+1-s)!
((k-1)n+t-1)!= \\
&& ((k-2)n+t-1)! \sum_{r=0}^s (-1)^{s-r}
\binom{(k-2)n+t-1 +r}{r} \binom{kn+t+1}{s-r} \times \nonumber \\
&& \ \ \ \ \ \ \ \ \ \ (r+(k-2)n+t-1)_{n+1} (r+(k-2)n+t-1)_{n+1}. \nonumber
\end{eqnarray}

Next we shall consider $V_{n,k}^{X,Y}$ and $A_{n,k}^{X,Y}$
where $X = Y= i+k\N$ for $0 \leq i \leq k-1$. In this case,
it is not longer the case that $A_{m,s}^{X,X} = V_{m,s}^{X,X}$ so
we will handle the cases of $A_{m,s}^{X,X}$ and $V_{m,s}^{X,X}$ separately.

First we shall consider  $V_{n,s}^{X,Y}$.  Note that if
$A = i +k\N$ and $B= i +k\N \cup i+1 +k\N$, then for
$m+1 = kn+i \in X \cap Y$, then
$x_n = |X \cap [m]| = n = y_m = |Y \cap [m]|$ and
$b_m = |B\cap [m]| = 2n$. Thus it follows from
Theorem \ref{D=V} that $V_{n,s}^{X,X} = D_{n,s}^{A,B}$ for
all $n$ and $s$. We then have
two cases. \\
\ \\
{\bf Case I.} $m = kn+t$ where $0 \leq t < i$. In this case, $x_m
=|X \cap [m]| =n$ and $x_n^c = |[m]-X| =(k-1)n+t$. Then it follows
from Theorem \ref{formula} that
\begin{equation}\label{1ap}
V_{m,s}^{X,X} = \binom{n}{s}^2 \binom{(k-1)n+t}{n-s} s! (n-s)! ((k-1)n+t)!.
\end{equation}

On the other hand, we can obtain a formula for
$V_{m,s}^{X,X}= D_{m,s}^{A,B}$ from equation (\ref{ap2}).
In this case $|A_m^c| = kn+t-n = (k-1)n+t$.
For any $0 \leq a \leq n$,
\begin{eqnarray}\label{2ap}
\alpha_{A,m,ak+i} +\beta_{B,m,ak+i} &=&  kn+t-1  - (n -(a+1)) - 2a
\nonumber \\
&=& (k-1)n+t-a.
\end{eqnarray}
Thus
\begin{eqnarray*}
&&\prod_{x \in A_m} (1+r +\alpha_{A,m,x} +\beta_{B,m,x}) = \\
&& \prod_{a=0}^{n-1} (1+r +\alpha_{A,m,ak+i} +\beta_{B,m,ak+i}) =\\
&& \prod_{a=0}^{n-1} (1+r +(k-1)n+t-a)  =\\
&&(1+r+(k-1)n+t)\downarrow_n
\end{eqnarray*}
where $(a)\downarrow_n$ is defined by $(a)\downarrow_0 =1$ and
$(a)\downarrow_n=a(a-1) \cdots (a-n+1)$ for $n \geq 1$. Thus it
follows that
\begin{eqnarray}
&&\binom{n}{s}^2 \binom{(k-1)n+t}{n-s} s! (n-s)!
((k-1)n+t)! = \\
&& ((k-1)n+t)! \sum_{r=0}^s (-1)^{s-r}
\binom{(k-1)n+t+r}{r} \binom{kn+t+1}{s-r} (1+r+(k-1)n+t)\downarrow_{n}. \nonumber
\end{eqnarray}

{\bf Case II.} $m = kn+t$ where $i\leq t \leq k-1$. In this case,
$x_m =|X \cap [m]| =n+1$ and $x_n^c = |[m]-X| =(k-1)n+t-1$. Then it
follows from Theorem \ref{formula} that
\begin{equation}\label{12ap}
V_{m,s}^{X,X} = \binom{n+1}{s}^2 \binom{(k-1)n+t-1}{n+1-s} s! (n+1-s)!
((k-1)n+t-1)!.
\end{equation}

On the other hand, we can obtain a formula for
$V_{m,s}^{X,X}= D_{m,s}^{A,B}$ from equation (\ref{ap2}).
In this case $|A_m^c| = kn+t-(n +1)= (k-1)n+t-1$.
For any $0 \leq a \leq n$,
\begin{eqnarray}\label{22ap}
\alpha_{A,m,ak+i} +\beta_{B,m,ak+i} &=&  kn+t-1  - (n+1 -(a+1)) - (2a +1)
\nonumber \\
&=& (k-1)n+t-1-a.
\end{eqnarray}
Thus
\begin{eqnarray*}
&&\prod_{x \in A_m} (1+r +\alpha_{A,m,x} +\beta_{B,m,x}) = \\
&& \prod_{a=0}^{n} (1+r +\alpha_{A,m,ak+i} +\beta_{B,m,ak+i}) =\\
&& \prod_{a=0}^{n-1} (1+r +(k-1)n+t-1-a)  =\\
&&(r+(k-1)n+t)\downarrow_{n+1}.
\end{eqnarray*}
Thus it follows
that
\begin{eqnarray}
&&\binom{n+1}{s}^2 \binom{(k-1)n+t-1}{n+1-s} s! (n+1-s)!
((k-1)n+t-1)!= \\
&& ((k-1)n+t-1)! \sum_{r=0}^s (-1)^{s-r}
\binom{(k-1)n+t-1+r}{r} \binom{kn+t+1}{s-r} (r+(k-1)n+t)\downarrow_{n+1}. \nonumber
\end{eqnarray}

Next we consider the case of computing  $A_{n,s}^{X,Y}$ where $X = Y
= i+k\N$ where $k \geq 2$ and $0 \leq i \leq k-1$. Let $A = i+k\N$
and $B=i-1+k\N$.  Then it is easy to see that if $m+1 =kn+i \in X
\cap Y = X$, then $x_m =y_m = n$ and $b_m =2n+1$.  Thus it follows
from Theorem \ref{D=A} that $A_{n,k}^{X,Y} = D_{n,s}^{A,B}$ for all
$n$ and
$s$ in this case. We then have two cases.\\
\ \\
{\bf Case A.} $m = kn+t$ where $0 \leq t <i -1$. In this case, $x_m
= n$ and $x_m^c = (k-1)n+t$. Thus if follows from Theorem \ref{AXX}
that
\begin{equation}
A_{m,s}^{X,X} =n! ((k-1)n+t)! \binom{n-1}{s} \binom{(k-1)n+t+1}{n-s}.
\end{equation}

On the other hand, we can obtain a formula for
$A_{m,s}^{X,X}= D_{m,s}^{A,B}$ from equation (\ref{ap2}).
In this case $|A_m^c| = (k-1)n+t$.
For any $0 \leq a \leq n$,
\begin{eqnarray}\label{22ap2}
\alpha_{A,m,ak+i} +\beta_{B,m,ak+i} &=&  kn+t-1  - (n-(a+1)) - (2a+1)
\nonumber \\
&=& (k-1)n+t-1-a.
\end{eqnarray}
Thus
\begin{eqnarray*}
&&\prod_{x \in A_m} (1+r +\alpha_{A,m,x} +\beta_{B,m,x}) = \\
&& \prod_{a=0}^{n-1} (1+r +\alpha_{A,m,ak+i} +\beta_{B,m,ak+i}) =\\
&& \prod_{a=0}^{n-1} (1+r +(k-1)n+t-1-a)  =\\
&&(r+(k-1)n+t)\downarrow_{n}.
\end{eqnarray*}
Thus it follows
that
\begin{eqnarray}
&&n! ((k-1)n+t)! \binom{n-1}{s} \binom{(k-1)n+t+1}{n-s}= \\
&& ((k-1)n+t)! \sum_{r=0}^s (-1)^{s-r}
\binom{(k-1)n+t+r}{r} \binom{kn+t+1}{s-r} (r+(k-1)n+t)\downarrow_{n}.
\nonumber
\end{eqnarray}
\ \\
{\bf Case B.} $m = kn+t$ where $i \leq t <k-1$. In this case, $x_m =
n+1$ and $x_m^c = (k-1)n+t-1$. Thus if follows from Theorem
\ref{AXX} that
\begin{equation}
A_{m,s}^{X,X} =(n+1)! ((k-1)n+t-1)! \binom{n}{s} \binom{(k-1)n+t}{n+1-s}.
\end{equation}

On the other hand, we can obtain a formula for
$A_{m,s}^{X,X}= D_{m,s}^{A,B}$ from equation (\ref{ap2}).
In this case $|A_m^c| = (k-1)n+t-1$.
For any $0 \leq a \leq n$,
\begin{eqnarray}\label{22ap3}
\alpha_{A,m,ak+i} +\beta_{B,m,ak+i} &=&  kn+t-1  - (n+1-(a+1)) - (2a+1)
\nonumber \\
&=& (k-1)n+t-2-a.
\end{eqnarray}
Thus
\begin{eqnarray*}
&&\prod_{x \in A_m} (1+r +\alpha_{A,m,x} +\beta_{B,m,x}) = \\
&& \prod_{a=0}^{n+1} (1+r +\alpha_{A,m,ak+i} +\beta_{B,m,ak+i}) =\\
&& \prod_{a=0}^{n+1} (1+r +(k-1)n+t-2-a)  =\\
&&(r+(k-1)n+t-1)\downarrow_{n+1}.
\end{eqnarray*}
Thus it follows
that
\begin{small}
\begin{eqnarray}
&&n! ((k-1)n+t-1)! \binom{n+1}{s} \binom{(k-1)n+t}{n+1-s}= \\
&& ((k-1)n+t-1)! \sum_{r=0}^s (-1)^{s-r}
\binom{(k-1)n+t-1+r}{r} \binom{kn+t+1}{s-r} (r+(k-1)n+t-1)\downarrow_{n+1}.
\nonumber
\end{eqnarray}
\end{small}

\section{Direction for future research}\label{beyond}
In this section, we shall describe some problems for further
research that naturally arise from the work in this paper.

There are other statistics which are closely related to the statistics
that we consider in this paper. For example, suppose that
$X,Y \subseteq \N$ and
define
$$\gamma_{X,Y}(\sg) =|\{i\in X:
\sg_i\in X\}\cup \{i\in Y: \sg_i\in Y\}|.$$
Let
$\Gamma_{n,s}^{X,Y}=|\{\sg\in\S_n:\gamma_{X,Y}(\sg)=s\}|$. Then we have
the following theorem.
\begin{theorem} For any $X$
and $Y$ such that $X\cup Y=\N$ and $X\cap Y=\emptyset$, we have
$\Gamma_{n,s}^{X,Y}=0$ unless $s=2k+y_n-x_n$ for some $k$, in which
case $$\Gamma_{n,2k+y_n-x_n}^{X,Y}=(x_n)!(y_n)!{x_n \choose k}{y_n
\choose x_n-k}.$$
\end{theorem}
\begin{proof}
Suppose we pick $k$ positions in $X_n$ to contain
elements in $X_n$ in  ${x_n \choose k}$ ways and we pick $x_n-k$
positions $Y_n$ to put the other elements of $X_n$ in ${y_n \choose
x_n-k}$ ways. The remaining positions in the permutation must be filled
with elements of $Y_n$. Next arrange elements in $X_n$ in $(x_n)!$
ways and we arrange elements in $Y_n$ in $(y_n)!$.
Clearly the number of permutations $\sg$ that can be constructed
in this way is $(x_n)!(y_n)!{x_n \choose k}{y_n
\choose x_n-k}$. Note that our construction
forces $y_n -(x_n -k)$ elements of $Y_n$ to be in positions
in $Y_n$ so that for any $\sg$ constructed in this way
$\gamma_{X,Y}(\sg) = 2k+y_n-x_n$.
\end{proof}

Note that in the special case where $X = \E$ and $Y = \O$, we have
that $\Gamma_{2s,2n}^{\E,\O} = (n!)^2 \binom{n}{s}^2$ and
$\Gamma_{2s+1,2n+1}^{\E,\O} = n!(n+1)! \binom{n}{s}\binom{n+1}{s+1}$
which agrees with other formulas in our table in these special cases.
However, for general $X$ and $Y$, we get quite different
recursions. For example, suppose that
$n+1 \in X \cap Y$, $\sg \in S_n$, $i \in X-Y$, and $\sg_i \in
Y-X$. Then it is easy to see that
$$\gamma_{X,Y}(I_{n+1}^{(i)}(\sg)) = \gamma_{X,Y}(\sg) +2$$
so that the value of $\gamma_{X,Y}$ can jump by 2 with a single
insertion. This type of phenomenon does not happen with any of
the other statistics studied in this paper. Thus it would be
interesting to further study these types of statistics to see if
one can prove explicit formulas of the type given in
Theorem \ref{HR1} and \ref{HR2}.

Even though we found solutions to all the bijective questions related to
the objects in our table, in some cases, one should be able to modify
our bijections (find new ones) to preserve more than one statistic.

Recall that the statistic $S_{10}$ is the number of odd
descent-tops, and $S_{12}$ is the number of (odd,odd) pairs. In this
section, we will use the following statistics as well.

\begin{itemize}
\item $S_{17}$ --- the maximal subsequence of the form $1 2\ldots i$ in a
permutation. E.g., $S_{17}(34152)=2$ while the increasing
permutation of length $n$ gives the maximum value of $S_{17}$ in
$\S_n$. A modification of this statistic was studied by
Zeilberger~\cite{Z} in connection with
       {\em 2-stack sortable permutations}.
\item $T_1 = S_{10}$ but {\em not} $S_{12}$.
\item $T_2 = S_{12}$ but {\em not} $S_{10}$.
\item $T_3 = S_{10}$ {\em and} $S_{12}$.

\end{itemize}

Our first conjecture is the following joint equidistribution:

\begin{conjecture}\label{con1}
The following should be true: $(S_{10},S_{12},S_{17}) \sim
(S_{12},S_{10},S_{17})$.
\end{conjecture}

Notice, that Conjecture~\ref{con1} suggests existence of an
involution turning $S_{10}$ to $S_{12}$ and vice versa. This
conjecture can be refined as follows.

\begin{conjecture}\label{con2} $(T_1, T_2, T_3, S_{17}) \sim (T_2, T_1, T_3, S_{17})$.
That is, if the involution mentioned above exists, it is likely to
leave pairs that are both $S_{10}$ and $S_{12}$ untouched.
\end{conjecture}

Observe that to preserve statistic $S_{17}$ in
Conjectures~\ref{con1} and~\ref{con2}, we need to require the
increasing $n$-permutation to go to itself, and this is the only
thing we need to worry about in our recursive construction of the
bijection regarding $S_{17}$ as otherwise it is not changed and thus
preserved by induction no matter where we stick the largest element.
So, it seems like we should be able to have the increasing
permutation as a fixed point.

Here is how a proof of Conjecture~\ref{con2} could be arranged for
even $n$ assuming the rest is constructed by induction. For odd
$n$'s things seem to be much more complicated.

For $n=1$, 1 is mapped to 1. Suppose we have constructed a bijection
from $\S_{n-1}$ to $\S_{n-1}$ ($n$ is even) such that it sends $k$
(resp. $\ell$, $s$) occurrences of $T_1$ (resp. $T_2$, $T_3$) to $k$
(resp. $\ell$, $s$) occurrences of $T_2$ (resp. $T_1$, $T_3$).
Inserting $n$ in $T_1$ (resp. $T_2$, $T_3$) pair decreases the
number of $T_1$ (resp. $T_2$, $T_3$) by 1 keeping all other
statistics unchanged. Clearly, we can manage the corresponding
insertion on the other side that would decrease by 1 the number of
occurrences of the corresponding statistic. Inserting $n$ at any
other position does {\em not} change a thing in either side and can
be matched to each other. In particular, inserting $n$ at the end
corresponds to inserting $n$ at the end and this guarantees that the
statistic $S_{17}$ is preserved (either it is unchanged in both
cases, or assuming we deal with the increasing permutation $12\ldots
(n-1)$ going to itself, $S_{17}$ is increased by 1 in both cases).

\end{document}